\documentstyle[11pt,fullpage,doublespace,amsmath,amsthm,amsfonts]{article}
\newtheorem{lemma}{Lemma}[section]
\newtheorem{proposition}[lemma]{Proposition}
\newtheorem{remark}[lemma]{Remark}
\newtheorem{theorem}[lemma]{Theorem}
\newtheorem{definition}[lemma]{Definition}

\newtheorem{example}[lemma]{Example}

\theoremstyle{remark}

\newcommand{\Id}{{{\mathchoice {\rm 1\mskip-4mu l} {\rm 1\mskip-4mu l}
      {\rm 1\mskip-4.5mu l} {\rm 1\mskip-5mu l}}}}

\begin{document}

\title {A Comparison of Hofer's Metrics on Hamiltonian 
Diffeomorphisms and Lagrangian Submanifolds}

\author{
Yaron Ostrover\thanks{This paper is a part of the author's 
Ph.D. 
thesis, being carried out under the supervision of Prof. Leonid Polterovich, at 
Tel-Aviv university.} 
\\
School of Mathematical Sciences \\
Tel Aviv University
69978 Tel Aviv, Israel \\
\tt yaronost@post.tau.ac.il}

\maketitle

\begin{abstract}
We compare Hofer's geometries on two spaces associated with
a closed symplectic manifold $(M,\omega)$. 
The first space is the group of Hamiltonian diffeomorphisms.
The second space ${\cal L}$ consists of all Lagrangian
submanifolds of $M \times M$ which are exact Lagrangian 
isotopic to the diagonal.
We show that in the case of a closed symplectic manifold with $\pi_2(M) = 
0$,
the canonical embedding of Ham$(M)$ into ${\cal L}, \   f 
\mapsto graph(f)$
is not an isometric embedding, although it preserves Hofer's 
length of smooth paths.
\end{abstract}

\section{Introduction and Main Results}

In this paper we compare Hofer's geometries on two remarkable
spaces associated with a closed symplectic manifold $(M,\omega)$.
The first space Ham$(M,\omega)$
is the group of Hamiltonian diffeomorphisms. 
The second consists of all Lagrangian submanifolds
of $(M \times M,- \omega \oplus \omega)$ which are exact Lagrangian
isotopic to the diagonal $\triangle \subset M \times M$.
Let us denote this second space by ${\cal L}$.
The canonical embedding $$j  : \text {Ham}(M,\omega) \rightarrow {\cal L},
 \ \ \  f \mapsto {\text graph}(f)$$
preserves Hofer's length of smooth paths.
Thus, it naturally follows to ask 
whether $j$ is an isometric embedding with
respect to
 Hofer's distance.
Here, we provide a negative answer to this question for the case of
a closed symplectic manifold with $\pi_2(M) =0$.
In fact, our main result shows that the image of
 Ham$(M,\omega)$ inside ${\cal L}$ is ``strongly distorted'' 
(see Theorem \ref{main_result} below).

Let us proceed with precise formulations. Given a path $\alpha = \{ 
f_t \}, \ t
\in [0, 1]$ of Hamiltonian diffeomorphisms of $(M,\omega)$, define
its Hofer's length (see [H]) as
$$ {\rm {length}(\alpha) } =
  \int_0^1 \Bigl\{ \max_{x \in M} F(x,t) -
\min_{x \in M} F(x,t) \Bigr\} \ dt$$
where $F(x,t)$ is the Hamiltonian function generating $\{ f_t \}$.
For two Hamiltonian diffeomorphisms $\phi$ and
$\psi$, define the Hofer distance  
$ d(\phi,\psi) = {\rm inf} \ {\rm length} (\alpha)$
where the infimum is taken over all smooth paths 
$\alpha$ connecting $\phi$ and 
$\psi$.
For further discussion see e.g. [LM1],[MS], and [P1].

Hofer's metric can be defined in a more general context of Lagrangian
submanifolds (see [C]).
Let $(P,\sigma)$ be a closed symplectic manifold, and let $\triangle
\subset  P$ be a closed Lagrangian submanifold. Consider a smooth family
$\alpha = \{ L_t \}, \ t \in [0,1]$ of Lagrangian submanifolds, such that
each $L_t$ is diffeomorphic to $\triangle$.
We call $\alpha$ an {\it exact path} connecting $L_0$ and $L_1$, if there
exists a smooth map $\Psi : \triangle \times [0,1] \rightarrow P$ such that 
for every $t$, $\Psi(\triangle \times \{ t \}) = L_t$, and in addition $\Psi^* \sigma = 
dH_t \wedge dt$ for some smooth function $H: \triangle \times [0,1]
\rightarrow {\mathbb R}$.
The Hofer length of an exact path is defined by
$$ {\rm {length}(\alpha)} = \int_0^1 \Bigl\{ \max_{x \in \triangle} H(x,t) -
\min_{x \in \triangle} H(x,t) \Bigr \} \ dt. $$
It is easy to check that the above notion of length is well-defined. 
Denote by ${\cal L}(P,\triangle)$ the space of all Lagrangian submanifolds
of $P$ which can be connected to $\triangle$ by an exact path.
For two Lagrangian submanifolds $L_1$ and $L_2$ in ${\cal 
L}(P,\triangle)$, define the Hofer distance $\rho$ on ${\cal L}(P,\triangle)$ as
follows: $\rho(L_1,L_2) = {\rm inf} \  {\rm length} (\alpha)$, where the infimum is
taken over all exact paths on ${\cal L}(P,\triangle)$ that connect 
$L_1$ and $L_2$.

\medskip In what follows we choose $P = M \times M$, $\sigma = - \omega \oplus
\omega$ and take $\triangle$ to be the diagonal of $M \times M$.
We abbreviate ${\cal L} = {\cal L}(P,\triangle)$ as in the beginning of the
paper. 
Based on a result by Banyaga [B], it can be shown that every smooth
path on ${\cal L}(P,\triangle)$ is necessarily exact.
Our main result is the following:

\begin{theorem} 
Let $(M,\omega)$ be a closed symplectic manifold with $\pi_2(M) =0$.
Then there exist a family $\{ \varphi_t \},\  t \in [0,\infty)$ in
Ham$(M,\omega)$ and a
constant
$c$ such that:

\begin{enumerate}
\item $ d(\Id,\varphi_t) \rightarrow \infty \ as \ t \rightarrow \infty$.
\item $ \rho( {\rm {graph}} (\Id), {\rm {graph}} (\varphi_t)) = c$.
\end{enumerate}
\label{main_result}
\end{theorem}

In fact, we construct the above family $\{ \varphi_t \}$ 
explicitly:
\begin {example}
{\rm
Consider an open set $B \subset M$. Suppose that there exists a 
Hamiltonian
diffeomorphism $h$ such that $h(B) \cap \text{Closure}\;(B) = 
\emptyset$.
By perturbing $h$ slightly, we
may assume that all the fixed points of $h$ are
non-degenerate.
Let $F(x,t)$, where $x \in M$, $t \in [0,1]$ be a Hamiltonian function such that
$F(x,t) = c_0 < 0$ for all $x \in M \setminus B, \ t \in 
[0,1]$.
Assume that $F(t,x)$ is normalized such that for every $t$,
$\int_MF(t,\cdot)\omega^n=0$.
We define the family $\{ \varphi_t \},\
t \in [0,\infty)$ by $\varphi_t =
h f_t$, where $ \{f_t
\}$ is the Hamiltonian flow generated by $F(t,x)$.
As we'll see below, the family $\{ \varphi_t \}$ satisfies the
 requirements of Theorem 1.1.}
\end{example}

Theorem \ref{main_result} has some corollaries:
\begin{enumerate}
\item The embedding of Ham$(M,\omega)$ in ${\cal L}$ is not isometric, rather, the image of 
Ham$(M,\omega)$ in ${\cal L}$ is highly
distorted. The minimal path between two graphs of Hamiltonian
diffeomorphisms in ${\cal L}$, might pass through exact Lagrangian
submanifolds which are not the graphs of any Hamiltonian
diffeomorphisms. Compare with the situation described in [M], where
it was proven that 
in the case of a compact manifold, the space of 
Hamiltonian deformations of
the zero section in the cotangent bundle is locally
flat in the Hofer metric.
\item 
The group of Hamiltonian diffeomorphisms of a closed 
symplectic manifold with $\pi_2(M)=0$ has an infinite diameter with
respect to Hofer's metric. 
\item Hofer's metric $d$ on Ham$(M,\omega)$ {\it does not} coincide with
the Viterbo-type metric on Ham$(M,\omega)$ defined by Schwarz in [S]. 
\end{enumerate}
As a by-product of our method we prove the following result (see Section 3 below):

\begin{theorem}
Let $(M,\omega)$ be a closed symplectic manifold with
$\pi_2(M)=0$.
Then there exists an element $\varphi$ in $( {\rm Ham}(M,\omega),d)$
which cannot be joined to the identity by a minimal geodesic.
\end{theorem}

The first example of this kind was established by 
Lalonde and McDuff [LM2] for 
the case of $S^2$.
\\

{\bf Acknowledgment.} I would like to express my deep gratitude to
my supervisor, Professor Leonid Polterovich, for his encouragement
and for many hours of extremely useful conversations.
I would also like to thank Paul Biran and Felix Schlenk for  
many fruitful discussions. 

\section{Proof of The Main Theorem}

In this section we prove Theorem \ref{main_result}.
Throughout this section  
let $(M,\omega)$ be a closed symplectic manifold with $\pi_2(M)=0$.
Let $\{ \varphi_t \}, \ t\in [0,\infty)$ the family of Hamiltonian diffeomorphisms defined 
in Example 1.2.
We begin with the following lemma which states that Hamiltonian
diffeomorphisms act as isometries on the space $({\cal L},\rho)$. 
The proof of the lemma follows immediately from the definitions.

\begin{lemma}
Let $\Gamma : \triangle \times [0,1] \rightarrow M \times M $ be an exact
Lagrangian isotopy in $\cal L$ and let $\Phi : M \times M \rightarrow M \times M$
be a Hamiltonian diffeomorphism. Then
$$ {\rm {length}} \{ \Gamma \} = {\rm {length}} \{\Phi 
\circ \Gamma \}.$$
In particular, $\rho(L_1,L_2) = \rho(\Phi(L_1),\Phi(L_2))$ for 
every $L_1,  L_2 \in {\cal L}$.
\end{lemma}
Next, consider the following exact isotopy of the Lagrangian embeddings 
$ \Psi : \triangle \times [0,\infty ) \rightarrow M \times M, \ 
\Psi(x,t) = (x,\varphi_t(x)).$
We denote by $L_t = \Psi(\triangle \times \{ t \})$ the graph of 
$\varphi_t = h f_t $ in $M \times M$.
The following proposition will be proved in Section 5 below.

\begin{proposition}
For every $t \in [0,\infty)$ there exists a Hamiltonian isotopy $\{ 
\Phi_s \}, \ s \in [0,t]$ of 
$M \times M$, such that 
$\Phi_s(L_0) = L_s$ and such that for every $s$, $\Phi_s(\triangle) =
\triangle$.
\end{proposition}
Hence, it follows from Proposition 2.2 and Lemma 2.1, that the family
$\{ \varphi_t \}, \ t \in [0,\infty )$ satisfies the second 
conclusion of Theorem 1.1 with constant $c = \rho(\triangle,L_0)$.
\\ \\
Let us now verify the first statement of Theorem 1.1.
For this purpose 
we will use a theorem by Schwarz [S] stated below. 
First, recall the definitions of the
action functional and the action spectrum.
Consider a closed symplectic manifold $(M,\omega)$ with $\pi_2(M) =0$. 
Let $\{ f_t \}$ be a Hamiltonian path generated by a Hamiltonian function
$F : [0,1]  \times M \rightarrow {\mathbb R}$.
We denote by ${\rm Fix}^{\circ}(f_1)$ the set of fixed points, $x$, of the 
time-1-map $f_1$ whose orbits $\gamma = \{ f_t(x) \} , \ t \in [0,1]$ 
are contractible.
For $x \in {\rm Fix}^{\circ}(f_1)$, take any 2-disc $\Sigma \subset M$ with
$\partial \Sigma = \gamma$, and define the symplectic action functional by
$$ {\cal A}(F,x) = \int_{\Sigma} \omega - \int_0^1 F(t,f_t(x))dt.$$
The assumption $\pi_2(M) = 0$ ensures that the integral $\int_{\Sigma} \omega$
does not depend on the choice of $\Sigma$.

\begin{remark}
{\rm
In the case of a closed symplectic manifold with $\pi_2(M) = 0$, a result
by Schwarz [S], implies that
for a Hamiltonian path $\{ f_t \}$ with $f_1 \neq \Id$ there exist two 
fixed points $x$, $y \in {\rm Fix}^{\circ}(f_1)$ with
${\cal A}(F,x) \neq {\cal A}(F,y)$. Moreover, 
the action functional does not depend on the choice of the Hamiltonian
path generating $f_1$.
Therefore, we can speak about the action of a fixed point of a 
Hamiltonian diffeomorphism, regardless of the Hamiltonian function
used to define it.}
\end{remark}

\begin{definition}
For each $f$ in $ {\rm Ham}(M,\omega)$ we define the action spectrum
$$ \Sigma_{f} = \{ {\cal A}(f,x) \ | \ x \in
{\rm Fix}^{\circ}(f) \} \subset {\mathbb R} . $$
\end{definition}
The action spectrum $\Sigma_{f}$ is a compact subset of $\mathbb R$ 
(see e.g. [S],[HZ]).
\begin{theorem}{{ \rm [S]}}.
Let $(M,\omega)$ be a closed symplectic manifold with $\pi_2(M) = 0$. Then, for
every  $f$ in {\rm Ham}$(M,\omega)$
$$ d(\Id,f) \geq \min \Sigma_{f}. $$
\end{theorem}

Next, consider the family $\{ \varphi_t \} = \{ hf_t \}, \ t \in [0,\infty)$.
Note that  
${\rm Fix}^{\circ}(h) = {\rm Fix}^{\circ}(\varphi_t)$ for every $t$.
The following proposition shows that the action spectrum
of $ \varphi_t $ is a linear translation of the action
spectrum of $h$. Its proof is carried out in Section 4.

\begin{proposition}
For every $t \in [0,\infty)$, and for every  
fixed point $z\in {\rm Fix^{\circ}}(\varphi_t) = {\rm 
Fix^{\circ}}(h)$,
$$ {\cal A}({\varphi_t},z) = {\cal A}(h,z) - tc_{0}$$
where $c_0$ is the negative (constant) value that $F$ attains on $M
\setminus B$ (see Example 1.2).  
\end{proposition}

We are now in a position to complete the proof of Theorem \ref{main_result}. Indeed, the action spectrum is a compact subset of $\mathbb R$, hence its minimum is finite.
By proposition 2.6 the minimum of $\Sigma_{\varphi_t}$ tends to infinity as 
$t \rightarrow \infty$. Thus,
$$ d(\Id,\varphi_t) \rightarrow \infty \ as \ t \rightarrow \infty $$
as follows from Theorem 2.5. This completes the proof of Theorem \ref{main_result}. 
\hfill $\square$

\section{Geodesics in Ham$(M,\omega)$ and Proof of Theorem 1.3}

In this section we describe our result about geodesics in the group
of Hamiltonian diffeomorphisms endowed with the Hofer metric $d$.
We refer the reader to [BP], [LM1], [LM2], and [P2] for further details 
on this subject.
\\ \\
Let $\gamma = \{ \phi_t\}, \ t \in [0,1]$ be a smooth regular path in
Ham$(M,\omega)$,
i.e. $\frac{d}{dt} \phi_t \neq 0$ for every $t\in [0,1]$. 
The path $\gamma$ is called a
{\em minimal geodesic} if it minimizes the distance between its
end-points:
$${\text {length}} (\gamma)=d(\phi_0,\phi_1).$$
The graph of a Hamiltonian path $\gamma = \{ \phi_t\}$ is the
family of embedded images of $M$ in $M \times M$ defined by
the map $ \Gamma : M \times [0,1] \rightarrow M \times M, \ (x,t) \mapsto 
(x,\phi_t(x))$.
Next, consider the family $\{ \varphi_t \}, \ t \in [0,\infty)$
that was constructed in Example 1.2.
We will show that there exists no
minimal geodesic joining the identity and $ \varphi_{t_0}$, for some
$t_{0}$.

\begin{proof}[Proof of Theorem 1.3]
Assume (by contradiction) that for every $t$, there exists 
a minimal geodesic in Ham$(M,\omega)$ joining the
identity with $ \varphi_t $.
Fix $t_0 \in [0,\infty)$. There 
exists a Hamiltonian path $\alpha = \{ f_s
\}, \ s \in [0,1]$ in
Ham$(M,\omega)$ such that
$$ d_{t_0} := d(\Id,\varphi_{t_0}) = {\rm
{length}}(\alpha).$$
Expressed in Lagrangian submanifolds terms, 
$ \Psi =\{ {\rm {graph}}(f_s) \}, \ s \in [0,1]$ is
an exact path in $M \times M$ joining the
diagonal with ${\rm {graph}}(\varphi_{t_0})$.
By Proposition 2.2, there exists a  Hamiltonian isotopy 
$\Phi$ such that for every $t$,
$\Phi_t{({\rm {graph}} (\varphi_{t_0})) = {\rm {graph}} (\varphi_{t})}$, and 
$\Phi_t(\triangle) = \triangle$.
We will choose $t_1$ to be sufficiently close to
$t_0$ so as to ensure that $\{ \Phi_{t_1}{({\rm {graph}}(f_s))} \}, \ s \in [0,1]$ is the 
graph of some
Hamiltonian path $\gamma$ in Ham$(M,\omega)$.
We claim the following
$$ d_{t_1} \leq {\rm {length}}(\gamma) 
= {\rm {length}} \{ {{\rm {graph}}(\gamma)}\} = 
 {\rm {length}} \{ {{\rm {graph}}(\alpha)} \} 
= {\rm {length}} \{ \alpha \}  =  d_{t_0}.  $$
Indeed, a straightforward computation yields that 
the embedding $f \mapsto graph(f)$ preserves Hofer's length, and from 
Lemma 2.1, ${\rm {length}} \{ {{\rm {graph}}(\alpha)} \} = 
{\rm {length}} \{ {{\rm {graph}}(\gamma)}\} $.
We have shown that for every $t_0$ there exists 
$\varepsilon > 0$ such that if $|t - t_0 | \leq
\varepsilon$ then $d_t \leq d_{t_0}$.
Since $d_t$ is a continuous
function, we conclude that $d_t$ is a constant function. 
On the other hand, by Theorem \ref{main_result}, $d_t = d(\Id,\varphi_t) \rightarrow \infty$ as $t \rightarrow \infty$.  
Hence there is a contradiction.
\end{proof}

\section{Proof of Proposition 2.6}


We investigate the expression ${\cal
A}({\varphi_t},z)$ for some fixed $t$. 
Since the action functional does not depend on the choice of the
Hamiltonian path generating the time-1-map (see Remark 2.3), 
we consider the following path generating $\varphi_t$.
$$ \gamma(s) =    \left\{ \begin{array}{ll}
        f_{2st} \ \ \ \ \ \  \ \ \ , & \mbox{ \ $s \in$ [0,{$\frac 1 2$}]  
} 
\\
        h_{2s-1} f_t \ \ \ \ ,  & \mbox{ \ $s \in$ ( {$\frac 1 2$},1].}
\end{array} \right.$$
Note that since $h(B) \cap B = \emptyset$ and $f_t$ is supported
in $B$, then for $z \in {\text Fix^{\circ}}(\varphi_t) = {\text 
Fix^{\circ}}(h)$ 
the path $\{ \gamma_s(z)\}, s \in [0,1]$ 
coincides with the path $\{ h_s(z) \}, \ s \in [0,1]$.
Denote by $\alpha$ the loop $\{ \gamma_s(z) \}, \ s \in [0,1]$ 
and let $\Sigma$ be any 2-disc with $\partial \Sigma = \alpha$.
The details of the calculation of 
${\cal A}({\varphi_t},z)$ are as follows:
$$ {\cal A}({\varphi_t},z) = \int_{\Sigma} \omega -
\int_{0}^{1}
t F(s,z)ds - \int_{0}^{1} H(s,h_s(z))ds,$$  
where $F$ and $H$ are the Hamiltonian functions generating $\{ h_t \}$ and 
$\{ f_t \}$
respectively.
Recall that by definition, $F$ is equal to a constant $c_0$ 
in $M \setminus B$.
This implies that  
$$ {\cal A}({\varphi_t},z) = \int_{\Sigma} \omega -
\int_{0}^{1} H(s,h_s(z))ds -tc_0. $$
The right hand side is exactly  ${\cal A}(h,z)-tc_0$.
Hence, the proof is complete. 



\section{Extending the Hamiltonian Isotopy}

In this section we prove Proposition 2.2. 
Let us first recall some relevant notations.
Let $\{ \varphi_t \}, \ t\in [0,\infty)$ the family of Hamiltonian 
diffeomorphisms defined in Example 1.2.
Consider the following exact isotopy of Lagrangian embeddings 
$ \Psi : \triangle \times [0,\infty ) \rightarrow M \times M, \
\Psi(x,t) = (x,\varphi_t(x)).$
We denote by $L_t = \Psi(\triangle \times \{ t \})$ the graph of
$\varphi_t = h f_t $ in $M \times M$, and by $\triangle$ the diagonal in 
$M \times M$.
It follows from the construction of the family $\{ \varphi_t
\}$, that for every $t$,  
${\rm Fix}(\varphi_t) = {\rm Fix}(h)$. Hence, 
 $L_t$ intersects the diagonal at the same set of points for
every $t$. 
Moreover, we assumed that all the fixed points of $h$ are
non-degenerate, therefore for every $t$, $L_t$ transversely intersect the diagonal.
In order to prove Proposition 2.2, we first need the following lemma.

\begin{lemma}
Let $x,y \in Fix^{\circ}(\varphi_t) = Fix^{\circ}(h)$,
i.e., intersection points of the family $\{ L_t \} $ and the diagonal in 
$M \times M$.
Take a smooth curve $\alpha : [0,1] \rightarrow M$
with $\alpha(0) = x$ and $\alpha(1) = y$ and
let $\Sigma : [0,1] \times [0,1] \rightarrow M, \ \Sigma(t,s) =
\varphi_t(\alpha(s))$ be a 2-disc such that
$\partial \Sigma_t = \varphi_t \alpha - h \alpha$.
Then the symplectic area of
$\Sigma_t = \Sigma(t,\cdot)$ vanishes for all $t$. 
\end{lemma}

\begin{proof}
By a direct computation of the symplectic area of
$\Sigma_t$, we obtain that
$$ \int_{\Sigma_t} \omega = \int_{ [0,t] \times [0,1]} \Sigma_t^* \omega =
- \int_0^t dt \int_0^1 dF_t( \frac {\partial} {\partial s} \varphi_t
\alpha(s)) \ ds
= \int_0^t {\widehat F_t}(\varphi_t(x))dt - \int_0^t {\widehat F_t} 
(\varphi_t(y))dt,$$
where ${\widehat F}$ is the Hamiltonian function generating the flow $\{ \varphi_t 
\}$. A straightforward computation shows that ${\widehat F}(t,x) = F(t,h^{-1}x)$, 
where $F$ is the Hamiltonian function generating the flow $\{ f_t \}$. 
Recall that by definition, $F(x,t)$ is equal to a constant $c_0$ outside  the ball $B$.
Moreover, since 
$x$, $y \in Fix^{\circ}(h)$ and $h(B) \cap B = \emptyset$, 
then $x,  y \notin$ B. 
Therefore, ${\widehat F_t}(\varphi_t(x)) =  {\widehat F_t}(\varphi_t(y)) = c_0$ for 
every $t$.
Thus, we conclude that for every $t$, the symplectic area 
of $\Sigma_t$ vanishes as required.
\end{proof}

\begin{proof}[Proof of Proposition 2.2]
We shall proceed along the following lines.
By the Lagrangian tubular neighborhood theorem (see [W]),
there exists a symplectic identification
between a small tubular neighborhood $U_s$ of $L_s$ in $M \times M$ and
a tubular neighborhood $V_s$ of the zero section in the
cotangent bundle $T^*L_s$.
Moreover, it follows from a standard compactness argument that 
there exists $\delta_s = \delta(s,U_s) > 0$ such that 
$L_{s'} \subset U_s$ for every $s'$ with $|s' - s| \leq \delta_s$.
Next, denote $I_s = (s - \delta_s, s + \delta_s) \cap [0,t]$, and
consider an open cover of the interval $[0,t]$ by the family $\{I_s \}$, that is
$ [0,t] = \bigcup_{s \in [0,t]} I_s$.
By compactness we can choose a finite number of points
$S = \{ s_1 < \ldots < s_n \}$ such that 
$[0,t] = \bigcup_{i = 1} ^ n I_{s_i}$.
Without loss of generality we may assume that $I_{s_j} \cap I_{s_{j+2}} = \emptyset$. 
Now, for every $s \in$ S,  
we will construct a Hamiltonian 
function ${\widetilde H}_s : U_s \rightarrow {\mathbb R}$ such 
that 
the corresponding 
Hamiltonian flow 
will shift $L_s$ toward $L_{s'}$ 
for $s' \in  I_s$,
and will leave the diagonal invariant.
Next, by smoothly patching together those Hamiltonian flows 
on the intersections $U_{s_i} \cap U_{s_{i+1}}$, we will achieve 
the required Hamiltonian isotopy $\Phi$. 
\\
\\
We fix $s_0 \in$ S. 
Let $(p,q)$ be  
 canonical local coordinates on $T^*L_{s_0}$ (where $q$ is the coordinate 
on $L_{s_0}$ and $p$ is the coordinate on the fiber). Moreover, we
fix a Riemannian metric on $ L_{s_0}$,
 and denote by
$\| \cdot \|_{s_0}$ the induced fiber norm on $T^*L_{s_0}$. 
Consider the aforementioned tubular neighborhood $ U_{s_0}  $ of 
$ L_{s_0} $ in $M \times M$.
For every $x \in L_{s_0} \cap \triangle$ denote by $\sigma_{s_0}(x)$
the component of the intersection of $U_{s_0}$ and 
$\triangle$ containing the point $x$.
Note that we may choose $U_{s_0}$ small enough such that the sets 
$\{ \sigma_{s_0}(x) \}, \ {x \in L_{s_0} \cap \triangle}$, are mutually disjoint.   
In what follows we shall denote the image of $\sigma_{s_0}(x)$ under the above 
identification between $U_{s_0}$ and $V_{s_0}$, by $\sigma_{s_0}(x)$ as well.
\\
\\
We first claim that there exists a Hamiltonian symplectomorphism 
${\widetilde \varphi} : V_{s_0} \rightarrow V_{s_0}$ which for every intersection point $x \in  
L_{s_0} \cap \triangle$
sends $\sigma_{s_0}(x)$ to the fiber over $x$ and which leaves $L_{s_0}$
invariant. Indeed, 
since $L_{s_0}$ transversely intersects the diagonal, and since $\sigma_{s_0}(x)$ is a 
Lagrangian submanifold,
 $\sigma_{s_0}(x)$ is the graph of a closed 1-form of $p$-variable
i.e, $\sigma_{s_0}(x) = \{ (p,\alpha(p)) \}$ where $\alpha(p)$ is 
locally defined near the intersection point $x$, 
and $\alpha(0) = 0$.
Define a family of local diffeomorphisms by
$\varphi_t (p,q) =  (p,q - t \alpha(p))$. Since the 1-form $\alpha(p)$ is
closed, $\{ \varphi_t \}$ is 
a Hamiltonian flow. 
Denote by $K(p,q)$ the Hamiltonian
function generating $\{ \varphi_t \}$. A simple computation shows that 
$K(p,q) =   
- \int 
\alpha(p) dp$. Hence $K(p,q)$ is independent on the $q$-variable 
i.e, 
$K(p,q) = K(p)$. Furthermore, we may 
assume that $K(0) = 0$.
 Next, we cut off the Hamiltonian function $K(p)$ outside 
a neighborhood of the
intersection point $x$.
Let $\beta(r)$ be a smooth cut-off function that vanishes for $r \geq 2 
\varepsilon$
and equal to 1 when $r \leq \varepsilon$, for sufficiently small 
$\varepsilon$. 
Define  $${\widetilde K}(p,q) = \beta (\| p \|) \cdot \beta(\| q \|) \cdot K(p).$$
A straightforward computation shows that, 
$ {\frac {\partial {\widetilde K}} {\partial 
q}}
(0,\cdot) =
 {\frac {\partial {\widetilde K}} {\partial p}} 
(0,\cdot) = 0$.
Hence
the time-1-map of the Hamiltonian flow corresponding to 
${\widetilde K(p,q)}$ is the required symplectomorphism.   
Therefore, we now can assume that 
$\sigma_{s_0}(x)$ coincide with the fiber over the point $x$.
\\
\\
Next, since $\Psi$ is an exact Lagrangian isotopy, we have that for 
every $s \in I_{s_0}$,
$L_{s}$ is a graph of an exact 1-form $dG_s$ in 
the symplectic tubular neighborhood $V_{s_0}$ of $L_{s_0}$. 
Hence, in the above local coordinates $(p,q)$ on $T^* L_{s_0}$, 
$L_s$ takes the form $L_s = 
(dG_s(q),q)$. Moreover, note that $dG_s(0) = 0$.
\\
\\
Define
%
$$ {\widetilde H}_{s_0}(p,q) = \beta(\|p\|)  \cdot G_s(q). $$ 
Consider the Hamiltonian vector field corresponding to ${\widetilde H_{s_0}}$,
$$ \widetilde {\xi} = \ \left \{ \begin{array}{ll}   
             \dot{p} = - \frac{\partial \widetilde {H}}{\partial q} = -
\beta(\|p\|)  \cdot \frac {\partial {G_s}(q)}{\partial q}
\\            \dot{q} =  \frac{\partial \widetilde {H}}{\partial p} =
  \frac{\partial } {\partial p} \{ \beta (\| p \|) \}            
\cdot {G_s}(q)      
       \end{array} \right .$$
It follows that for every $s \in I_{s_0}$ such that $L_s \subset 
 \left  \{ (p,q)  \ | \ \|  p \| <  {\varepsilon}  \right \}$, the
Hamiltonian flow is given by $$ (p,q) \rightarrow \left ( p + \frac 
{\partial 
G_s(q)} {\partial q}   ,q \right )$$
Hence, locally, the Hamiltonian flow shift $L_{s_0}$ toward 
$L_s$ as required.
It remains to prove that $\widetilde {\xi}$ vanishes on the diagonal.
First, since $dG_s(0) = 0$, it follows that $\dot {p} = 0$.
Next, consider $x$ and $y$, two intersection points of the family 
$\{ L_s \}$ and the diagonal. 
It follows from Lemma 5.1 that 
the symplectic area between  
$L_{s_0}$ and $L_s$ in $V_{s_0}$ vanishes for every $s \in I_{s_0}$.
Hence, by the same argument as in Lemma 5.1, for every such $s$ we have 
$$ 0 = \int_{\Sigma_s} \omega = \int_{ [0,s] \times [0,1]} \Sigma_s^* 
\omega 
= \int_0^s \Bigl ( G_s(x) -  G_s(y) \Bigr )  ds$$
Thus, we get that $G_s(x) - G_s(y) = 0$.  
Note that by changing the functions $\{ G_s \}$ by a summand 
depending
only on s, we can assume that for every $s$, ${G_s}$ vanishes on $L_s
\cap \triangle$.
It now easily follows that $\widetilde {\xi}_{ |_\triangle}=0$.
%
%
%
Therefore, we have that the diagonal is invariant 
under 
the Hamiltonian flow. 
Finally, by (smoothly) patching together all the 
Hamiltonian flows corresponding to the Hamiltonian functions 
${\widetilde H_{s_i}}$, for $i=1, \ldots ,n$, 
 we conclude that  
 there exists a Hamiltonian isotopy $\Phi$ 
such that $\Phi_s(L_0) = L_s$ 
and $\Phi_s(\triangle) = \triangle$. 
This completes the proof of the proposition.
\end{proof}


\centerline{\bf References}

\smallskip

\begin{description}

\item {[B]} Banyaga, A. (1978). {\it Sur la structure du groupe des 
diff\'eomorphisms qui pr\'eservent une forme symplectique.} Comm.
Math. Helv. 53, 174-227.

\item{[BP]} Bialy, M. and Polterovich, L. {\it
Geodesics of Hofer's metric on the group of Hamiltonian
diffeomorphisms}, Duke Math. J., 76 (1994), 273--292.

\item{[C]} Chekanov, Yu. {\it Invariant Finsler metrics on the space
of Lagrangian embeddings}, Math. Z., 234 (2000),
605--619.   

\item{[H]} Hofer, H. {\it On the topological properties of symplectic
maps.} Proceedings of the Royal Society of Edinburgh, 115 (1990), 25-38.

\item{[HZ]} Hofer, H. and Zehnder, E. {\it Symplectic invariants and
Hamiltonian dynamics}, Birkhauser Advanced Texts, Birkhauser Verlag, 1994.

\item{[LM1]} Lalonde, F. and McDuff, D. {\it The geometry of symplectic
energy}, Ann. of Math 141 (1995), 349-371.

\item{[LM2]} Lalonde, F. and McDuff, D. {\it Hofer's $L^\infty$ - 
geometry:
energy and stability of Hamiltonian flows}, parts 1 and 2,
Invent. Math. 122 (1995), 1-33 and 35-69.

\item{[M]} Milinkovi\'c , D. {\it Geodesics on the space of Lagrangian
submanifolds in cotangent Bundles}, Proc. Amer. Math. Soc. 129 (2001), no 
6, 1843--1851. 

\item{[MS]} McDuff, D. and Salamon, D. {\it Introduction to
Symplectic Topology}, 2nd edition, Oxford University Press, Oxford,
England
(1998).   

\item{[P1]} Polterovich, L. (1993). {\it Symplectic displacement energy
for Lagrangian submanifolds.} Ergodic Theory and Dynamical Systems, 13,
357-67.

\item{[P2]} Polterovich, L. {\it The Geometry of the group of 
Symplectic
Diffeomorphisms}, Lectures in Math, ETH, Birkhauser (2001).

\item{[P3]} Polterovich, L. {\it Growth of maps, distortion in groups 
and symplectic geometry}, Preprint Math.DS/0111050.

\item{[S]} Schwarz, M. {\it On the action spectrum for closed
symplectically aspherical manifolds}, {\it Pac. Journ. Math} {\bf 193}
(2000),
419--461.

\item{[W]} Weinstein, A. (1971). {\it Symplectic manifolds and their Lagrangian submanifolds.} Advances in Mathematics, 6, 329-46.

\end{description}
\end{document}